\documentclass{amsart}
\title{Notes on two methods for direct construction of probabilistic LFSR sequences of third order}
\author{Lhoussain El Fadil and Danilo Gligoroski} 
\thanks{Partially supported by ERCIM}
\newfont{\titre}{cmbx12 at 16pt}
\newtheorem{proposition}{Proposition}

\newcommand{\z}{{Z}\!\!\!{Z}}

\newcommand{\N}{{I}\!\!\!{N}}
\newcommand{\C}{{I}\!\!\!{C}}
\newcommand{\q}{{I}\!\!\!{Q}}
\newcommand{\cqfd}
{%
\mbox{}%
\nolinebreak%
\hfill%
\rule{2mm}{2mm}%
\medbreak%
\par%
}

\newcommand{\biindice}[3]%
{

\begin{array}[t]{c}
#1\\
{\scriptstyle #2}\\
{\scriptstyle #3}
\end{array}

}

\date{}
\begin{document}
 \maketitle
 \section*{Introduction}
In \cite{Tm}, the authors give two public key encryptions based on
third order linear sequences modulo $n^2$, where $n=pq$ is an RSA
integer. In their scheme (3), there are two mistakes in the
decryption procedure:
\begin{enumerate}
\item
The owner of the private key does not know the value of $m$ such
that $C_1=s_m(a,b)$ and $C_2=s_{-m}(a,b)$, and thus he/she can not
compute $L(C_1,C_2)$.
\item
If $s_{\lambda}(a,b)= 3$ modulo $n^2$, then  $L(a,b)$ is not
invertible modulo $n$. It follows that the owner of the private key can
not decrypt the cipher $C=(C_1,C_2)$ since he/she can not compute
$\frac{L(C_1,C_2)}{L(a,b)}$.
\end{enumerate}
In this short note, in order to decrypt the ciphertext
$C=(C_1,C_2)$, another map $L$ similar to that given in \cite{Tm} is
constructed. More precisely, if $L(a,b)$ is not invertible modulo
$n$, we describe a method how to choose $(a,b)$ such that  $L(a,b)$
is invertible modulo $n$ and how to compute $L(C_1,C_2)$.

 \section{Third order Linear sequences}
Let $p$ be an odd prime integer, $(a,b)\in\z^2$ and  $s(a,b)$ be a
the third linear order sequence  defined by
 $s_{k+3}(a,b) =as_{k+2}(a,b)-bs_{k+1}(a,b)+s_k(a,b)$ ( $(a,b)$ is
 called the generator of $s(a,b)$ and $k$ is the exponent).
 Let $\alpha_1, \alpha_2$ and $\alpha_3$ be the complex roots of $f(X)$.
 Then there exists $(a_1,a_2,a_3)\in\q^3$ such that for every
 $k\in z$,  $s_k(a,b)=a_1\alpha_1^k+a_2\alpha_2^k+a_3\alpha_3^k$.
Note that the tuple $(a_1,a_2,a_3)$ depends on the choice of
$s_{2}(a,b), s_0(a,b)$ and $s_1(a,b)$.
 For $(a_1,a_2,a_3)\in\z^3$ such that $a_1=a_2=a_3=1\,[p]$ (modulo $p$), we have
$s_0(a,b)=3$, $s_1(a,b)=a$ and $s_{-1}(a,b)=b$ modulo $p$.

In the following, we assume that $a_1=a_2=a_3=1\,[p]$. Denote   $T_{K/\q}$ and $N_{K/\q}$ the trace and norm maps of $K$.  Then $\overline {s_k(a,b)}= \overline{T_{K/\q}(\alpha_1^k)} \,[p]$ and
for every $i$, $\overline {N_{K/\q}(\alpha_i)}={(\overline \alpha_i)}^{p^2+p+1}=1 \,[p]$. Thus, $p^2+p+1$ is a period of $s(a,b)$ modulo $p$.\\

 The following cryptographic  applications of LFSR sequences are listed in
\cite{SL, G}. We present them without proof.
\begin{enumerate}
 \item
   For every  $k\in\z$, let $f_k(X)=X^3- s_k(a,b)X^2+s_{-k}(a,b)X-1\,[p]$.  Then $f_k(X)=(X-\alpha_1^k)(X-\alpha_2^k)(X-\alpha_3^k)\,[p]$.
\item
In particular, for every  $k$ and $e$,  $s_e(s_k(a,b), s_{-k}(a,b))= s_{ke}(a,b)\,[p]$.
\end{enumerate}
\section{Clarification remarks considering \cite{Tm}}
 In this section, in order to decrypt a  cipher $C=(C_1,C_2)$
 as given in  \cite[Sec.3 and 4]{Tm}, let $n=pq$ be an RSA, $(a,b)\in \N^2$
 such that  $f(X)=X^3-aX^2+bX-1$ is irreducible modulo $p$ (resp. modulo $q$).
 Let $s_k(a,b)$   be the third order linear sequence  modulo $n^2$ generated by $(a,b)$
 such that $s_0(a,b)=3$, $s_1(a,b)=a$ and $s_{-1}(a,b)=b$ modulo $n$.\\

 In order to have $\frac{s_{\lambda}(a,b)- 3}{n} \,[n]$ invertible,
 if   $s_{\lambda}(a,b)\neq 3\,[n^2]$,  then we will keep   $s_0(a,b)=3$, $s_1(a,b)=a$ and $s_{-1}(a,b)=b$ modulo $n^2$.
 If $s_{\lambda}(a,b)= 3\,[n^2]$, then  let $A=a+n$, $\beta_1$, $\beta_2$ and $\beta_3$ the roots of $g(X)=X^3-AX^2+bX-1$. Let $s(A,b)$ be the characteristic sequence  generated by
$A$ and $b$  modulo $n^2$ :
$\left \{ \begin{array}{cccc}
    s_{k+3}(A,b)&=&As_{k+2}(A,b)-bs_{k+1}(A,b)+s_{k}(A,b)\,[n^2]\\
     s_{0}(A,b)=3, & s_{1}(A,b)=A, &s_{-1}(A,b)=b\,[n^2]
\end{array}
\right.$

Since $\bar{f}(X)=\bar{g}(X) \,[n]$, then up to a permutation for every $1 \le i\le 3$, there exists $t_i\in\C$ such that $\beta_i=\alpha_i+nt_i$. Thus, for every integer $k$,
 $s_{k}(A,b)=\sum_{i=1}^3\beta_i^k=\sum_{i=1}^3(\alpha_i+nt_i)^k=\sum_{i=1}^3(\alpha_i^k+nkt_i)= s_{k}(a,b)+nk(t_1+t_2+t_3)$ modulo $[n^2]$. For $k=1$, we have $t_1+t_2+t_3=1\,[n]$. Thus,
$s_{\lambda}(A,b)=s_{\lambda}(a,b)+n\lambda \,[n^2]$, and then  $\frac{s_{\lambda}(A,b)- 3}{n}=\lambda \,[n]$ is invertible.
Finally, without loss of generality, up to replace $a$  by $a+n$, we can assume that  $\frac{s_{\lambda}(a,b)-3}{n}$
 is invertible modulo $n$, where  $\lambda$ is the least common multiple of $(p^2+p+1,q^2+q+1)$.

\begin{proposition}
\begin{enumerate}
 \item
   For every  $k\in\z$, let $f_k(X)=X^3- s_k(a,b)X^2+s_{-k}(a,b)X-1\,[n^2]$.  Then $f_k(X)=(X-\alpha_1^k)(X-\alpha_2^k)(X-\alpha_3^k)\,[n^2]$.
\item
In particular, for every  $k$ and $e$,  $s_e(s_k(a,b), s_{-k}(a,b))= s_{ke}(a,b)\,[n^2]$.
\end{enumerate}
\end{proposition}
{\bf Proof.}
 Since $s_{k}(a,b)=\sum_{i=1}^3\alpha_i^k\,[n^2]$, $\alpha_1^k\alpha_2^k\alpha_3^{k}=1\,[n^2]$ and $\alpha_1^k\alpha_2^k+\alpha_1^k\alpha_3^k+\alpha_2^k\alpha_3^k=\alpha_1^{-k}+\alpha_2^{-k}+\alpha_3^{-k}\,[n^2]$.\cqfd

Let $\Gamma=\{(x,y)\in \z^2,\, s_{\lambda}(x,y)=3\,[n]\}$ and  $L:\,\Gamma \longrightarrow \frac{\z}{n\z}$ be the map defined
 by $L(x,y)=\frac{s_{\lambda}(x,y)-3}{n}\, [n]$. Since   $ s_{\lambda}(x,y)=3\,[n]$, then $L$ is well defined.

 \begin{proposition}
For every integer $k$, $\frac{L(s_{k}(a,b), s_{-k}(a,b))}{L(a,b)}=k\,[n]$.
\end{proposition}
{\bf Proof.}
First, $L(a,b)=\frac{s_{\lambda}(a,b)-3}{n}$ is invertible modulo $n$.
 Let $\Gamma_{1}^i=\{x\in \frac{\z[\alpha_i]}{n^2\z[\alpha_i]},\, x=1\,[n]\}$ and  ${L}^i:\,\Gamma_1^i \longrightarrow \frac{\z[\alpha_i]}{n\z[\alpha_i]}$
 be the map defined by ${L}^i(x)=\frac{x-1}{n}\, [n]$. Then for every $(x,y)\in(\Gamma_1^i)^2$, ${L}^i(xy)=\frac{xy-1}{n}=\frac{x(y-1)+(x-1)}{n}=x\frac{(y-1)}{n}+\frac{(x-1)}{n}=x{L}^i(y)+{L}^i(x)\, [n]$.   Since $x=1\, [n]$,
${L}^i(xy)={L}^i(x)+{L}^i(y)\, [n]$.\\
 Let $\lambda=k_p(p^2+p+1)$ and $1\le i\le 3$. Since $\overline {N_{K/\q}(\alpha_i)}=\alpha_i^{p^2+p+1}\, [p]$,  $\alpha_i^{\lambda}=(\overline {N_{K/\q}(\alpha_i)})^{k_p}=1\, [p]$ (resp. modulo $q$). It follows that $\alpha_i^{\lambda}=1\, [n]$,  $\alpha_i^{\lambda}\in\Gamma_{1}^i$ and  ${L}^i(\alpha_i^{k\lambda})=k{L}^i(\alpha_i^{\lambda})\,[n]$.
  Therefore, $s_{k\lambda}(a,b)-3=(\alpha_1^{k\lambda}-1)+(\alpha_2^{k\lambda}-1)+(\alpha_3^{k\lambda}-1)\,[n]$, $L(s_{k}(a,b), s_{-k}(a,b))=\frac{s_{k\lambda}(a,b)-3}{n}= \sum_{i=1}^3{ L}^i(\alpha_i^{k\lambda})=k\sum_{i=1}^3{L}^i(\alpha_i^{\lambda})\,[n]$ and $L(a,b)=\frac{s_{\lambda}(a,b)-3}{n}=\sum_{i=1}^3{L}^i(\alpha_i^{\lambda})\,[n]$. As $L(a,b)$ is invertible modulo $n$, $\frac{L(s_{k}(a,b))}{L(a,b)}=k\,[n]$.\cqfd

\subsection{ The deterministic version}
{\bf Algorithm of encryption and decryption in Scheme 3 of
\cite{Tm}}
\begin{enumerate}
\item
Public parameters: $(n, a,b)$
\item
Private parameters: $(p,q)$
\item
Encryption: For  a message $0\le m < n$,  Bob calculates the
ciphertext block $C=(c_1,c_2)$ such that  $c_1=s_m(a,b),
c_2=s_{-m}(a,b)\,[n^2]$.
\item
Decryption: For  a given ciphertext block $c$, Alice can decrypt it by calculating $\frac{L(c_1,c_2)}{L(a,b)}\,[n]$. \\
 \end{enumerate}

Indeed, since $c=(c_1,c_2)$ is a ciphertext, let $0\le m < n$ such
that $c_1=s_{m}(a,b)$ and $c_2=s_{-m}(a,b)$. Then
$s_{\lambda}(c_1,c_2)=s_{\lambda}(s_m(a,b),s_{-m}(a,b))=s_{m\lambda}(a,b)=3$
modulo $n$. Thus,  $L(c_1,c_2)$ is well defined and
$\frac{L(c_1,c_2)}{L(a,b)}=\frac{L(s_{m}(a,b)}{L(a,b)}=m\,[n]$.\cqfd

 \subsection{ The probabilistic version}
{\bf Algorithm of encryption and decryption in Scheme 3 of
\cite{Tm}}
\begin{enumerate}
\item
Public parameters: $(n, a,b)$
\item
Private parameters: $(p,q)$
\item
Encryption: For  a message $ 0\le m < n$,  Bob selects a random
integer $r$ and calculates the ciphertext block $C=(c_1,c_2)$ such
that $c_1=s_{rn+m}(a,b), c_2=s_{-(rn+m)}(a,b)\,[n^2]$.
\item
Decryption: For  a given ciphertext block $c$, Alice can decrypt it by calculating $\frac{L(c_1,c_2)}{L(a,b)}\,[n]$.\\
 \end{enumerate}

 \end{document}